\documentclass[12pt,final]{amsart}

\newcommand{\R}{\mathbb R}
\newcommand{\Z}{\mathbb Z}

\newtheorem{thm}{Theorem}

\title{The constant of interpolation}
\author{Artur Nicolau}
\address{Dept. Matem\`atiques, Universitat Aut\`onoma de Barcelona, 08193
Bellaterra, Spain}
\email{artur@mat.uab.es}
\author{Joaquim Ortega-Cerd\`a}
\address{Dept.\ Matem\`atica Aplicada i An\`alisi,
Universitat  de Barcelona, Gran Via 585, 08071 Barcelona, Spain}
\email{quim@mat.ub.es}
\thanks{The authors are supported by the European Commission Research Training
Network HPRN-CT-2000-00116. The first two authors are supported by 
DGICYT grants:
BFM2002-00571, BFM2002-04072-C02-01 and by the CIRIT: 2001SGR00172,
2001SGR00431. The third author is
partly supported by a grant from the Research Council of Norway.}
\author{Kristian Seip}
\address{Dept. of Mathematical Sciences, Norwegian University of Science
and  Technology, N--7491 Trondheim, Norway}
\email{seip@math.ntnu.no}
\date{January 28, 2003}

\begin{document}
\begin{abstract}
We prove that a suitably adjusted version of Peter Jones' formula for
interpolation in $H^\infty$ gives a sharp upper bound for what is known
as the constant of interpolation. We show how this leads to precise and
computable numerical bounds for this constant.
\end{abstract}
\maketitle
With each finite or infinite sequence $Z=(z_j)$ ($j=1,2,...$) of
distinct points $z_j=x_j+iy_j$ in the upper
half-plane of the complex plane, we associate a number
$M(Z)\in \R^+\cup\{+\infty\}$ which we call the constant of interpolation.
We may define it in two equivalent ways. The first is related to Carleson's
interpolation theorem for $H^\infty$ \cite{Carleson58}.
We say that $Z$ is an interpolating
sequence if the interpolation problem
\begin{equation} \label{problem} f(z_j)=w_j, \ j=1,2,...
\end{equation}
has a solution $f\in H^\infty$
for each bounded sequence $(w_j)$
of complex numbers. Using the open mapping theorem,
we find that if $Z$ is an interpolating sequence, then we can always solve
\eqref{problem} with a function $f$ such that
\[
\|f\|_\infty \le C \|(w_j)\|_\infty
\]
for some $C<\infty$ depending only on $Z$.
The constant of interpolation $M(Z)$ is declared to be the smallest
such $C$. We set $M(Z)=+\infty$ if $Z$ is not an interpolating sequence.

By a classical theorem of Pick (see \cite[p.~2]{Garnett81}),  we may
alternatively define $M(Z)$   as follows. Let $M_n(Z)$ be the smallest number
$C$ such that the  matrices
\[
\left(\frac{1-\overline{w_j}w_k}
              {z_j-\overline{z_k}}\right)_{j,k=1,2,\ldots,n}
\]
are positive semi-definite whenever $\|(w_j)\|_\infty\le 1/C$. The
constant of interpolation is then
$M(Z)=M_n(Z)$ if $Z$ is a finite sequence consisting of $n$ points and
$M(Z)=\lim_{n\to\infty} M_n(Z)$ if $Z$ is infinite.
We will make no use of this definition,
but have stated it to make the reader aware of the relevance of $M(Z)$ for
the classical Nevanlinna-Pick problem.

Carleson's interpolation theorem \cite{Carleson58} states that
$Z$ is an interpolating sequence (or alternatively $M(Z)<\infty$)
if and only if
\[
\delta(Z)=\inf_{j\neq k }  \prod_{k\ne j}
\left|\frac{z_j-z_k} {z_j-\overline{z_k}}\right|>0.
\]
Clearly, an interpolating sequence satisfies the Blaschke condition.
We let $B$ be the associated Blaschke product and set
\[
B_j(z)=\frac{z-\overline{z_j}}{z-z_j} B(z),
\]
so that we may write $\delta(Z)=\inf_j|B_j(z_j)|$.

An interesting result related to Carleson's theorem is that if
$M(Z)<\infty$, then the interpolation may be obtained by means of a
linear operator. In fact, P.~Beurling \cite{Carleson63} proved that
there exist $f_j\in H^\infty$ with $f_j(z_j)=1$ and $f_j(z_k)=0$ if
$k\ne j$, such that
\[
M(Z)=\sup_z \sum_j |f_j(z)|.
\]
The functions $f_j$ have the form
\[
f_j(z) = \frac{B_j
(z)}{B_j (z_j)} \Bigl(\frac{2iy_j }{z-
\overline{z_j}}\Bigr)^2 \, \frac{G(z_j)}{G(z)},
\]
where $G$ is a bounded analytic function solving a certain
nonlinear extremal problem. Unfortunately, $G$ is not given
explicitly, and it seems very difficult to get much further.
The problem of finding $G$ can be seen as a version of
the Nevanlinna-Pick interpolation problem, where one is interested
in computing $M(Z)$ and finding solutions of minimal
norm. There are classical results of
R.~Nevanlinna describing these solutions, but they are very implicit and
give little help in concrete situations.
It is therefore of interest to find more explicit solution operators,
along with good estimates for $M(Z)$.

A remarkably simple formula was found by P.~Jones \cite{Jones83}.
He showed that the series
\[
f(z)=\sum_j w_j
\frac{B_j(z)}{B_j(z_j)} \Bigl(\frac{ 2i y_j}
{z- \overline{z_j}}\Bigr)^2 \exp \Bigl( -a i \sum_{y_k\le y_j}
\bigl( \frac{y_k}{z - \overline{z_k}} -
\frac{y_k}{z_j - \overline{z_k}} \bigr)\Bigr)
\]
defines a function $f\in H^\infty$ such that
$f(z_j)= w_j$ with $\|f\| \le  C \| (w_j)\|_\infty$.
Here $a$ can be chosen freely and $C$ is a constant depending on
$a$ and the
sequence $Z$. The purpose of this note is to show that
this explicit operator, conveniently adjusted, is close to
optimal. By considering a certain extreme configuration of points, we
are in fact able to prove that it yields a sharp upper bound for
$M(Z)$. As a result, $M(Z)$ may be
bounded from above and below by fairly explicit numerical constants.

We begin by showing how to
``optimize'' Jones' formula. Take an analytic function
$g$ such that $g(i)=1$. We need $|g|$ to have a harmonic majorant, so we
require $(z+i)^{-2}g(z)\in H^1$ \cite[p.~60]{Garnett81}. Let $u$ denote the
least harmonic majorant of $|g|$ and set
\[
g_j(z)=g((z-x_j)/y_j), \ \ u_j(z)=u((z-x_j)/y_j).
\]
We assume further that $g$ is such that
\[
U_k(z)=\sum_{y_j\le
y_k}\frac{u_j(z)}{|B_j(z_j)|}
\]
defines a harmonic function; let
$V_k(z)$ be a harmonic conjugate of $U_k$, and set $G_k=U_k+iV_k$.
This leads us to the following interpolation formula:
\[
f(z)=\sum_j w_j \frac{B_j(z)}{B_j(z_j)}
g_j(z)\exp\Bigl(- a (G_j(z)-G_j(z_j))\Bigr)
\]
with $a$ some constant which may be chosen freely. Clearly,
$f(z_j)=w_j$. We define
\[
c_J(Z,g)=\sup_j U_j(z_j),
\]
so that for arbitrary $z$ we get the estimate
\[
|f(z)|\le \|(w_j)\|_\infty \frac{\exp(a c_J(Z,g))}{a}
\sum_j \frac{a |g_j(z)|}{|B_j(z_j)|}\exp(-a U_j(z)).
\]
Replacing $|g_j|$ by $u_j$, we find that the latter sum is a
lower Riemann sum for the integral
\[
\int_0^\infty e^{-t}dt
\]
so that we arrive at the estimate
\[
|f(z)|\le \|(w_j)\|_\infty \frac{\exp(a c_J(Z,g))}{a}.
\]
We see that the optimal choice of $a$ is $1/c_J(Z,g)$, and this
leads us to the bound
\[
M(Z)\le e c_J(Z,g).
\]
We may finally minimize $c_J(Z,g)$ and define
\[
c_J(Z)=\inf_g c_J(Z,g)
\]
so that
\[
M(Z)\le e c_J(Z).
\]
We have then proved one part of the following theorem.

\begin{thm} For every sequence $Z$ in the upper half-plane,
\begin{equation}\label{star}
  M(Z)\le e c_J(Z).
\end{equation}
The inequality is best possible in the sense that the constant $e$
on the right side of \eqref{star} cannot be replaced by any smaller
number.
\end{thm}

We postpone for the moment the proof of the sharpness of
\eqref{star}; it will be established by means of an explicit example
at the end of this note.

It may be argued that finding the $g$ minimizing $c_J(Z,g)$ is not much
easier than solving for the function $G$ in P.~Beurling's formula.
However, we will now point out that $c_J(Z)$ relates nicely to more
computable characteristics.

An immediate observation is that if we choose $g(z)=-4/(z+i)^2$, then
$u(z)=4(y+1)/|z+i|^2$ so that $c_J(Z,g)$ becomes
\[
c_{HJ}(Z)=\sup_n \sum_{y_j\le
y_n} \frac{4y_j(y_j+y_n)}{|z_j-\overline{z_n}|^2}
\frac{1}{|B_j(z_j)|}.
\]
This choice of $g$ corresponds to the original
version of Jones' formula. (The letter `$H$' in $c_{HJ}(Z)$ stands for
Havin; see below.) For this characteristic we have the following result.

\begin{thm} For every sequence $Z$ in the upper half-plane,
\[
M(Z)\le k c_{HJ}(Z)
\]
for some universal constant $k$. The best possible $k$
lies in the interval $[\pi/\log4, e]=[2.2662...,2.7183...]$.
\end{thm}
We have already established the upper bound for $k$. The lower bound will
again follow from the example to be considered below.

Our third and final characteristic was introduced by V. Havin in the first
appendix of \cite{Koosis98}. We get it from the expression for
$M(Z)$ obtained from Carleson's duality argument
(see \cite[p.~135]{Garnett81}):
\[ \label{constint}
M(Z)=\sup\left\{4\pi\sum
\frac{y_j |h(z_j)|}{|B_j(z_j)|}: h\in H^1, \|h\|_1\le
1\right\}.
\]
If we choose $h(z)=\pi^{-1}y_k/(z-\overline{z_k})^2$, $k=1,2,...$,
we arrive at
\[
c_H(Z) = \sup_k \sum_j \frac {4y_ky_j}{|z_k-\bar z_j|^2}
\frac 1{|B_j(z_j)|}
\]
along with the estimate
\[
M(Z)\ge c_H(Z).
\]
Since clearly $c_{HJ}(Z)\le 2 c_{H}(Z)$,
we may summarize our findings as a chain
of inequalities:
\begin{equation}
\label{inequalities}
c_H(Z)\le M(Z) \le e c_J(Z)\le e c_{HJ} (Z)\le 2e c_H(Z).
\end{equation}

In \cite{Koosis98}, Havin proves that
\[
c_H(Z) \le M(Z) \le k c_{H}(Z),
\]
with $k$ a universal constant. To prove the right inequality, he
proceeds by duality and uses the invariant Blaschke
characterization of Carleson measures, which is closely related to
the original proof of Carleson. By computing both $c_H(Z)$ and $M(Z)$ when
$Z$ consists of two points, he also shows that the left inequality is
best possible. In fact, it may be checked that each of the
inequalities in our chain \eqref{inequalities} is sharp.

To interpret the ``geometric'' contents of our characteristics, it may be
useful to relate them to the condition
\begin{equation} \label{blaschke}
\sup_k \sum_j \frac{y_j y_k}{|z_k-\bar z_j|^2}<+\infty,
\end{equation}
which is called the invariant Blaschke condition
(see \cite[p.~239]{Garnett81}). We see that our three characteristics are
closely related to the supremum appearing in \eqref{blaschke}. It may
also be noted that
by the bound $M(Z)\le 2e c_H(Z)$ and a calculus argument
applied to the invariant Blaschke sum, we obtain
\[
M(Z)\le \frac{2e+4e\log (1/\delta(Z))}{\delta(Z)};
\]
see \cite[p.~268]{Koosis98}.

We finally turn to our example which proves the sharpness of \eqref{star}
and the lower bound for $k$ in Theorem 2. In what follows the
notation $a(\gamma)\sim b(\gamma)$ will mean that $a(\gamma)$ and
$b(\gamma)$ are asymptotically equal, i.e., $\lim_{\gamma\to +\infty}
a(\gamma)/b(\gamma)=1$.
\vspace{.2cm}

\textbf{An example}. Fix $\gamma>0$ and consider the Blaschke
product defined by
\[
B(z)=B (\gamma , z) = \prod_{k \le 0}
\frac{z-ie^{k/\gamma}}{z+ie^{k/\gamma}} \prod_{k > 0}
\frac{ie^{k/\gamma}-z}{z+ie^{k/\gamma}}.
\]
The signs have been chosen so that $iB'(i)>0$, which ensures
convergence of the product.
The sequence of zeros $Z_\gamma=(ie^{k/\gamma})_{k\in \Z}$ is
clearly an interpolating sequence with $M(Z_\gamma)$ blowing up
when $\gamma$ tends to $+\infty$. To obtain appropriate estimates for $B$,
we relate it to the function
\[
F(z) = 2 e^{-\frac{\pi^2\gamma}{2}}\sin (\pi\gamma \log (-iz)),
\]
where $\log(z)$ is the principal branch of the logarithm. Both
$B$ and $F$ are bounded functions, and they have the same zeros. The
quotient $F(z)/B(z)$ is an outer function with modulus
close to 1 when $\gamma$ is large. More precisely,
we have
\[
\sup_{x\in \R\setminus\{0\}}
\left|\log \frac{|F(x)|}{|B(x)|}\right|
\sim e^{-\pi^2\gamma},
\]
and therefore the same asymptotic relation holds
in the upper half-plane. The Blaschke product $B$ is highly
symmetric. It is real on the imaginary half-axis $i\R^+$  and
moreover $B(e^{1/\gamma} z)=-B(z)$. We check that on $i\R^+$
the modulus of $B$ peaks  at the points
$\{ie^{(k+1/2)/\gamma} : k\in\Z \}$.
Again comparing it to $F$, we check that \begin{equation}
\label{solve}
B(ie^{(k+1/2)/\gamma}) =
(-1)^k 2e^{-\frac{\pi^2\gamma}{2}} t_\gamma \ \ \text{with} \ \
t_\gamma \sim 1. \end{equation}

We will now obtain a lower estimate for $M(Z_\gamma)$ by finding a minimal
norm solution of
the interpolation problem
\[
f(ie^{k/\gamma})=(-1)^k,\ k\in \Z.
\]
By \eqref{solve}, the problem is solved by the function
\[
g(z)=c_\gamma B(e^{1/(2\gamma)} z),
\]
with $c_\gamma$ an appropriate constant satisfying
$c_\gamma \sim e^{\frac{\pi^2\gamma}{2}}/2$. This means that if
we can prove that $g$ is a minimal norm solution, then it follows
that
\begin{equation} \label{lower}
M(Z_\gamma)\ge \frac{t_\gamma}{2}e^{\frac{\pi^2\gamma}{2}} \ \ \text{with} \ \
t_\gamma \sim 1.
\end{equation}

We wish to prove that $g$ is a solution of minimal norm. To this end, observe
that an arbitrary minimal norm solution can expressed as
\[
f = g + h B
\]
with $h$ a bounded analytic function. We may assume that $f$ is real
on $i\R^+$ because by symmetry we may if necessary replace $f$ by
$(\overline{f(-\overline{z})}+f(z))/2$.
Thus $h$ is also real on $i\R^+$. We define
\[
h_m(z) = \frac 1m \sum_{k=0}^{m-1} h(e^{2k/\gamma}z),
\]
and choose a convergent subsequence $h_{m_k}(z)\to \tilde h(z)$
such that the limit function satisfies
$\tilde h(e^{2/\gamma} z)=\tilde h(z)$, and $\tilde
h(iy)\in \R$ for real $y$. Hence $\tilde f = g + \tilde h B$ is
also a minimal norm solution and $
\tilde f ( e^{2/\gamma} z)=f(z) $. Finally, note that
\[
\varphi(z)=\frac 12 (\tilde f(z) -
\tilde f(e^{1/\gamma} z)),
\]
is a minimal norm solution as well such that
\begin{equation} \label{period}
\varphi(e^{1/\gamma}z)=-\varphi(z).
\end{equation}

Assume now that $g$ is not a minimal norm solution. Then
$\|\varphi\|_\infty < \|g\|_\infty$. Between the points $i$ and
$ie^{1/\gamma}$, $\varphi$  has a zero $i\delta$ because it is real on $i\R$.
Therefore, by the periodicity expressed by \eqref{period},
$\varphi$ has zeros at
$i\delta e^{k/\gamma}$,  $k\in \Z$. It follows that we may
factorize $\varphi$ as
\[
\varphi(z) = B(z / \delta) \varphi_0(z).
\]
We evaluate $\varphi$ at the point $i$ and get
\[
1=|\varphi(i)| = |B(i / \delta)| |\varphi_0(i)| \leq \frac1{c_\gamma}
\|\varphi_0 \|_\infty =  \frac {\|\varphi\|_\infty}{c_\gamma} =
\frac{ \|\varphi \|_\infty }{\|g \|_\infty}< 1,
\]
which is a contradiction. We conclude that $g$ has minimal norm so that
\eqref{lower} holds.

The next step is to compute $c_{J}(Z_\gamma)$. Since $B(e^{1/\gamma}z)=-B(z)$,
we have that
\[
|e^{k/\gamma} B'(ie^{k/\gamma})| =
|B'(i)|\]
for each integer $k$. Hence
\[
|B_k (ie^{k/\gamma})| = 2
e^{k/\gamma}|B'(ie^{k/\gamma})| = 2|B'(i)|.
\]
The derivative
$B'(i)$ can be estimated in terms of $F'(i)$, which gives us
\[
iB'(i) e^{\frac{\pi^2\gamma}{2}}/(2\pi\gamma)\to 1 \
\text{as} \  \gamma\to +\infty.
\]
Thus
\begin{equation}
\label{first}
c_{J}(Z_\gamma) \sim (4\pi \gamma)^{-1}
e^{\frac{\pi^2 \gamma}{2}} \inf_{g(i)=1} \sup_{k\in \Z} \sum_{y_j\le
y_k} u(i y_k/y_j)
\end{equation}
with $u$ denoting as before the least harmonic majorant of $|g|$.
Using the explicit expression for this majorant, we get
\[
\inf_{g(i)=1} \sup_{k\in \Z} \sum_{y_j\le
y_k} u(i y_k/y_j)
= \inf_{g(i)=1} \sum_{k \ge 0}\frac{1}{\pi}\int_{\R} \frac{1}{1+t^2}
\, |g(e^{k/\gamma}t)| \, dt.
\]
We interpret the sum on the right as a Riemann sum, so that
\[
\begin{split}
\sum_{y_j\le y_k} u(i y_k/y_j)
\sim &
\frac{\gamma}{\pi} \int_\R
\frac 1{1+t^2}\int_0^\infty |g(te^x)|\, dx\, dt
= \\
&\frac{\gamma}{\pi}\int_\R
\frac 1{1+t^2}\int_t^\infty \frac{|g(u)|}{u}\, du\, dt.
\end{split}
\]
Integrating by parts, we get
\[
\sum_{y_j\le
y_k} u(i y_k/y_j)\sim
  \frac{\gamma}{\pi} \int_{\R} \frac{\arctan t}{t}\, |g(t)|\, dt.
\]
We want to minimize the latter integral over all functions
$g$ such that $(z+i)^{-1}g\in H^1$ and $g(i)=1$.
This can be restated as an extremal problem in the weighted Hardy space
with norm
\[
\|h\|^2 = \int_\R |h(t)|^2 \frac{\arctan t}t\, dt.
\]
In turn, we can reduce this problem to one for the standard Hardy space
$H^2$, and we find that our original problem is solved by the function
\[
g_0(z) = \left(\frac {2i}{z+i}\right)^2\frac{\psi(z)}{\psi(i)},
\]
where $\psi(z)$ is the outer function whose modulus is
  $t / \arctan t$ on $\R$. Since
\[
\int_{\R} |g_0(t)| \frac{\arctan t}{t}\, dt =
\int_{\R} \frac 4{(t^2+1)} \frac 1{|\psi(i)|}=\frac{4\pi}{|\psi(i)|},
\]
we  get
\begin{equation} \label{almost}
c_{J}(Z_\gamma)\sim \frac{e^{\frac{\pi^2 \gamma}{2}}}{\pi |\psi(i)|}
\end{equation}
when plugging our extremal function $g_0$ into \eqref{first}.

We are left with the computation of $|\psi(i)|$. We first note that
\[
\psi(i) = i\exp\Bigl(-\frac i\pi \int_{\R}\,\Bigl(\frac 1{i-t}+
\frac{t}{t^2+1}\Bigr)
\log|\arctan(t)|\, dt \Bigr).
\]
Since
\[
|\arctan(t)|=\frac 12 \Bigl|\log\Bigl(\frac{1-it}{1+it}\Bigr)\Bigr|,
\]
the change of variables
\[
e^{i\theta}=\frac{1-it}{1+it}
\]
brings us to the explicit expression
\[ \psi(i) = 2i\exp\Bigl(-\frac{1}{2\pi}\int_{-\pi}^{\pi} \log|\log
e^{i\theta}|\, d\theta \Bigr) = \frac{2ie}{\pi}.
\]
Combining \eqref{lower} and \eqref{almost}, we conclude that
\[
c_J(Z_\gamma)\sim \frac{1}{2e}
e^{\frac{\pi^2 \gamma}{2}} \le \frac{t_\gamma}{e} M(Z_\gamma) \ \
\text{with} \ \ t_\gamma\sim 1,
\]
which proves the sharpness of \eqref{star} of Theorem 1.

The computation of $c_{HJ}(Z_\gamma)$ is straightforward. Indeed,
\[
c_{HJ}(Z_\gamma)\sim (4\pi \gamma)^{-1}
e^{\frac{\pi^2 \gamma}{2}} \sum_{k\ge 0}
\frac{4e^{-k/\gamma}}{(1+e^{-k/\gamma})}.
\]
The sum is again regarded as a Riemann sum, i.e.,
\[
  \sum_{k\ge 0}
\frac{e^{-k/\gamma}}{1+e^{-k/\gamma}}
\sim \gamma\,  \int_0^{\infty}
\frac{e^{-x}} {1+e^{-x}}\, dx= \gamma\, \log 2
\]
so that we arrive at the relation
\[
c_{HJ}(Z_\gamma)\sim \frac{\log 2}{\pi}
e^{\frac{\pi^2 \gamma}{2}} \le \frac{t_\gamma 2 \log 2}{\pi} M(Z_\gamma)  \ \
\text{with} \ \ t_\gamma\sim 1.
\]
This proves the lower bound for $k$ in Theorem 2.

\providecommand{\MR}{\relax\ifhmode\unskip\space\fi MR }
\providecommand{\MRhref}[2]{%
  \href{http://www.ams.org/mathscinet-getitem?mr=#1}{#2}
}


\begin{thebibliography}{Koo98}

\bibitem[Car58]{Carleson58}
Lennart Carleson, \emph{An interpolation problem for bounded analytic
   functions}, Amer. J. Math. \textbf{80} (1958), 921--930. \MR{22 \#8129}

\bibitem[Car63]{Carleson63}
\bysame, \emph{Interpolations by bounded analytic functions and the {C}orona
   problem}, Proc. Internat. Congr. Mathematicians (Stockholm, 1962), Inst.
   Mittag-Leffler, Djursholm, 1963, pp.~314--316. \MR{31 \#549}

\bibitem[Gar81]{Garnett81}
John~B. Garnett, \emph{Bounded analytic functions}, Pure and Applied
   Mathematics, vol.~96, Academic Press Inc. [Harcourt Brace Jovanovich
   Publishers], New York, 1981. \MR{83g:30037}

\bibitem[Jon83]{Jones83}
Peter~W. Jones, \emph{{$L\sp{\infty }$} estimates for the {$\bar \partial $}
   problem in a half-plane}, Acta Math. \textbf{150} (1983), no.~1-2, 137--152.
   \MR{84g:35135}

\bibitem[Koo98]{Koosis98}
Paul Koosis, \emph{Introduction to {$H\sb p$} spaces}, second ed., Cambridge
   Tracts in Mathematics, vol. 115, Cambridge University Press, Cambridge, 1998,
   With two appendices by V. P. Havin [Viktor Petrovich Khavin].
   \MR{2000b:30052}

\end{thebibliography}
\end{document}